\documentclass[12pt]{iopart}
\usepackage{iopams}
\usepackage{theorem}
\usepackage{algorithm}
\usepackage{graphicx}
\usepackage[usenames]{color}
\usepackage{multirow}
\usepackage{dsfont}

\newcommand{\eqref}[1]{\eref{#1}}

\newcommand{\R}{\mathbb{R}}
\newcommand{\N}{\mathbb{N}}

\newcommand{\ve}[1]{\boldsymbol{#1}}

\def\xk{\ve x^{(k)}}
\def\x{{\ve x}}
\def\y{{\ve y}}
\def\z{{\ve z}}

\newcommand{\argmin}[1]{\arg\min_{#1}}

\def\arg{{\rm arg}}
\def\dk{\ve d^{(k)}}
\def\hatx{\hat{\ve x}}
\def\xstar{\ve x^*}
\def\Dkk{{D_k^{-1}}}
\def\ak{\alpha_k}

\def\yk{\ve y^{(k)}}\def\yj{\ve y^{(j)}}

\def\lamk{\lambda^{(k)}}
\def\lamj{\lambda^{(j)}}
\def\xkk{\ve x^{(k+1)}}
\def\xj{\ve x^{(j)}}
\def\xjj{\ve x^{(j+1)}}

\def\endproof{\hfill$\square$\vspace{0.2cm}\\\noindent}

\newtheorem{lemma}{Lemma}[section]
\newtheorem{theorem}{Theorem}[section]
\newtheorem{proposition}{Proposition}[section]

\newtheorem{definition}{Definition}[section]

  {\end{list}}

\makeatletter

  \gdef\listctr{list\romannumeral\the\@listdepth}\expandafter
  
\makeatother

  {\end{list}}

\begin{document}
\title[New convergence results for the scaled gradient projection method]{New convergence results for the scaled gradient projection method}
\author{S Bonettini$^1$ and M Prato$^2$}
\address{$^1$ Dipartimento di Matematica e Informatica, Universit\`a di Ferrara, Via Saragat 1, 44122 Ferrara, Italy}
\address{$^2$ Dipartimento di Scienze Fisiche, Informatiche e Matematiche, Universit\`a di Modena e Reggio Emilia, Via Campi 213/b, 41125 Modena, Italy}
\eads{\mailto{silvia.bonettini@unife.it}, \mailto{marco.prato@unimore.it}}
\begin{abstract}
The aim of this paper is to deepen the convergence analysis of the scaled gradient projection (SGP) method, proposed by Bonettini et al. in a recent paper for constrained smooth optimization. The main feature of SGP is the presence of a variable scaling matrix multiplying the gradient, which may change at each iteration. In the last few years, an extensive numerical experimentation showed that SGP equipped with a suitable choice of the scaling matrix is a very effective tool for solving large scale variational problems arising in image and signal processing. In spite of the very reliable numerical results observed, only a weak, though very general, convergence theorem is provided, establishing that any limit point of the sequence generated by SGP is stationary. Here, under the only assumption that the objective function is convex and that a solution exists, we prove that the sequence generated by SGP converges to a minimum point, if the scaling matrices sequence satisfies a simple and implementable condition. Moreover, assuming that the gradient of the objective function is Lipschitz continuous, we are also able to prove the ${\mathcal O}(1/k)$ convergence rate with respect to the objective function values. Finally, we present the results of a numerical experience on some relevant image restoration problems, showing that the proposed scaling matrix selection rule performs well also from the computational point of view.
\end{abstract}
%
\section{Introduction}
Several inverse problems in applied sciences can be addressed by means of a constrained optimization problem
\begin{equation}\label{minf}
\min_{\ve x\in \Omega} f(\ve x),
\end{equation}
where $\Omega\subseteq\R^n$ is a closed and convex set and $f$ is a continuously differentiable function. First order methods are attractive for solving \eqref{minf} especially when $n$ is large and when the Hessian is not available or difficult to exploit. Indeed, the main strengths of these methods are, in general, the low memory requirement and the low computational cost per iteration. When the constraints set $\Omega$ has some special structure (e.g. box constraints, simplexes, balls), gradient projection (GP) methods have shown to be a valid tool to solve \eqref{minf} in a variety of frameworks, such as signal and image processing \cite{Bardsley2006,Bertero2008}, statistical inference \cite{Figuereido2007,Lin2007} and machine learning \cite{Serafini2005a,Serafini2005b,Zanni2006}.\\
The increasing popularity of first order methods gave rise in the recent literature to several studies aiming to devise suitable approaches for improving the convergence properties. In particular, we mention the extrapolation/inertial techniques \cite{Bertsekas2009,Beck2009b,Nesterov2005} and the variable metric approach \cite{Birgin2003,Bonettini2009,Combettes2014}. In the first case, an extrapolation step ensures the ${\mathcal O}(1/{k^2})$ convergence rate of the objective function values to the optimal one, where $k$ is the iteration index. In the latter case, an acceleration of the progress towards the solution is achieved by adopting a variable metric, which, at each iteration, could better capture the local features of the problem.\\
The focus of this work is on the scaled gradient projection (SGP) method \cite{Bonettini2009}, a variable metric algorithm based on a Armijo line--search.
The basic SGP iteration is given by
\begin{equation}\label{iter0}
\xkk = \xk + \lamk \dk = \xk + \lamk(\yk-\xk),
\end{equation}
where $\yk$ is the scaled Euclidean projection of $\xk-\alpha_kD_k\nabla f(\xk)$ onto $\Omega$, i.e.
\begin{equation}\label{y1}
\yk = \argmin{\x\in\Omega} \nabla f(\xk)^T(\x-\xk) + \frac{1}{\ak}(\x-\xk)^T\Dkk(\x-\xk),
\end{equation}
$\alpha_k > 0$ is the stepsize parameter, $D_k$ is a symmetric positive definite scaling matrix and $\lamk\in(0,1]$ is computed by means of a linesearch backtracking procedure to guarantee the sufficient decrease of the objective function. In this scheme, $\alpha_k$ and $D_k$ have to be considered as 'free' parameters which, when chosen in a clever way, can significantly improve the convergence behaviour of the algorithm (see e.g. \cite{Barzilai1988,Dai2006a,DeAsmundis2013,Fletcher2012}). In particular, the recent literature shows that a suitable combination of the stepsize parameter $\alpha_k$ and of the scaling matrix $D_k$ makes SGP a very effective tool in solving convex \cite{Bonettini2009,Zanella2009,Bonettini2010a,Prato2012} and nonconvex \cite{Bonettini2011,Bonettini2015,Bonettini2013b,Prato2013} problems arising in signal and image processing applications.\\
However, the convergence analysis of SGP available in the literature only establish that, when $\alpha_k$ and the eigenvalues of $D_k$ are bounded above and below away from zero, any limit point of the sequence $\{\xk\}_{k\in\N}$ is stationary for problem \eqref{minf}. This result has been proved in \cite[Theorem 2.1]{Bonettini2009} without any further assumption and it is mainly based on the properties of the Armijo linesearch.\\
In this paper we provide a new, stronger, convergence result for SGP when applied to convex problems, establishing the convergence of the sequence  $\{\xk\}_{k\in\N}$ to a solution of \eqref{minf}, provided that the eigenvalues of $D_k$ converge to one as $k$ diverges, at a certain rate.\\
Moreover, if we further assume that the gradient of $f$ is Lipschitz continuous, we provide a new ${\mathcal O}(1/k)$ complexity result on the objective function value.\\
We observe that the ${\mathcal O}(1/k)$ complexity result is worse than the one obtained for the inertial/extrapolation methods such as the celebrated FISTA \cite{Beck2009b}. However, we also show that the practical performances of SGP are comparable to FISTA on some significant image restoration problems. Our numerical experience also shows that the condition on the scaling matrix selection ensuring the theoretical convergence of the method is also useful from a computational point of view.\\
The paper is organized as follows: in section \ref{sec:line-search} we recall the basic properties of the Armijo linesearch procedure and of descent methods. The analysis of the scaled gradient projection method is performed in section \ref{SGPconv}, where the relationship between our approach and the related literature is also discussed, while section \ref{sec:num} is devoted to some illustrative numerical examples. Our conclusions are given in section \ref{sec:concl}.
\paragraph{Notation and basic definitions}\label{par:1}
In the following $\|\cdot\|$ indicates the $\ell_2$ norm of a vector while $\|\cdot\|_{D}$ denotes the norm induced by the symmetric positive definite matrix $D$, i.e. $\|\x\|_D = (\x^TD\x)^{\frac 1 2}$; $\R_{>0}^n$ and $\R_{\geq 0}^n$ denote the positive and non-negative orthants of $\R^n$, respectively; $\mu_{\min}(A)$, $\mu_{\max}(A)$ are the minimum and maximum eigenvalue of a square matrix $A$, respectively. The notation $A\succeq B$, where $A,B\in\R^{n\times n}$ are symmetric positive semidefinite matrices, indicates that $A-B$ is positive semidefinite. Given $\mu\geq 1$, we denote by ${\mathcal M}_\mu$ the set of the symmetric positive definite matrices with all eigenvalues contained in the interval $[\frac 1 \mu,\mu]$. For any $D\in {\mathcal M}_\mu$ we have that $D^{-1}$ also belongs to ${\mathcal M}_\mu$ and
\begin{equation}\label{ine_norm}
\frac 1\mu\|\x\|^2 \leq \|\x\|^2_{D}\leq \mu\|\x\|^2 \qquad \forall \x\in\R^n.
\end{equation}
We also recall the definitions of {\it stationary point} and {\it descent direction} for problem \eqref{minf} (see for example \cite{Bertsekas1999}).
\begin{definition}\label{Def_stationary}
A point $\x\in \Omega$ is a stationary point for problem \eqref{minf} if
\begin{equation*}
\nabla f(\ve x) ^T (\ve y - \ve x) \geq 0  \quad \forall \ve y\in \Omega.
\end{equation*}
\end{definition}
\begin{definition}
Let $\ve x$ be any point of the set $\Omega$.
\begin{itemize}
\item[(i)] A vector $\ve d\in \R^n$ is a feasible direction at $\ve x$ if $\ve x +\ve d\in \Omega$.
\item[(ii)] A vector $\ve d\in \R^n$ is a descent direction at $\ve x$ for problem \eqref{minf} if it is feasible and $\nabla f(\ve x)^T\ve d < 0$.
\end{itemize}
\end{definition}
Finally, we report the definitions of convex, globally and locally Lipschitz and level bounded function.
\begin{definition}
A continuously differentiable function $f:\Omega \rightarrow \R$ is said to be:
\begin{itemize}
\item[(i)] convex if
\begin{equation}\label{def:convex}
f(\ve y) \geq f(\ve x) + \nabla f(\ve x)^T(\ve y - \ve x) \quad \forall \ve x, \ve y \in \Omega;
\end{equation}
\item[(ii)] globally Lipschitz, if
\begin{equation}
\|f(\ve y) - f(\ve x)\| \leq L\|\ve y - \ve x\| \quad \forall \ve x, \ve y \in \Omega;
\end{equation}
\item[(iii)] locally Lipschitz, if for every compact set $K \subseteq \Omega$ there exists $L_K > 0$ such that
\begin{equation}
\|f(\ve y) - f(\ve x)\| \leq L_K\|\ve y - \ve x\| \quad \forall \ve x, \ve y \in K;
\end{equation}
\item[(iv)] level bounded, if the set $\Omega_\zeta=\{\x\in\Omega:f(\x)\leq \zeta\}$ is bounded for every $\zeta\in\R$.
\end{itemize}
\end{definition}
\section{General results about Armijo based gradient projection methods}\label{sec:line-search}
In this section, we recall the basic properties of the most popular linesearch procedure, the {\it Armijo linesearch}, given in Algorithm \ref{Algo2}. These results allow to prove a general convergence result which applies to any method where the objective function over two successive iterates decreases at least as it would decrease by applying the Armijo linesearch procedure along a suitable descent direction.\\

\begin{algorithm}[ht]\caption{Armijo linesearch (LS) algorithm}\label{Algo2}
Let $\{\xk\}_{k \in \N}$ be a sequence of points in $\Omega$ and $\{\dk\}_{k \in \N}$ a sequence of descent directions. Choose some $\delta, \beta \in(0,1)$ and compute $\lambda^{(k)}$ as follows:
\begin{itemize}
\item[1.] Set $\lambda^{(k)}=1$
\item[2.] \textsc{If}
\begin{equation}
f(\xk+\lambda^{(k)}\dk)\leq f(\xk)+\beta \lambda^{(k)}\nabla f(\xk)^T\dk
\label{Armijo}
\end{equation}
\textsc{Then} go to step 3\\
\textsc{Else} set $\lambda^{(k)}=\delta\lambda^{(k)}$ and go to step 2
\item[3.] \textsc{End}
\end{itemize}
\end{algorithm}

\noindent For the linesearch procedure based on the Armijo rule we recall the following basic theorem, which can be derived from known results \cite{Bertsekas1999,Grippo2000}.
\begin{proposition}\label{ProArm}
Let $\{\ve x^{(k)}\}_{k\in\N}$ be a sequence of points in $\Omega$. Assume that $\ve x^{(k)}$ converges to some $\bar{\ve x}\in\Omega$ and let $\{\ve d^{(k)}\}_{k\in\N}$ be a sequence of descent directions such that
\begin{equation}\label{A2}
\nabla f(\ve x^{(k)})^T \ve d^{(k)} < 0 \quad \forall k \in \N.
\end{equation}
Then the LS algorithm is well defined, i.e. for each $k\in\N$ it terminates in a finite number of steps. If, in addition, there exists a number $M>0$ such that $\|\ve d^{(k)}\|\leq M \ \forall k \in \N$ and
\begin{equation}\label{A3}
\displaystyle\lim_{k\rightarrow\infty}f(\xk)- f(\xk+\lambda^{(k)}\ve d^{(k)})=0,
\end{equation}
where $\lambda^{(k)}$ is computed with Algorithm \ref{Algo2}, then we have
\begin{equation*}
\lim_{k\rightarrow\infty} \nabla f(\ve x^{(k)})^T\ve d^{(k)} = 0.
\end{equation*}
\end{proposition}
It is worth stressing that the previous proposition applies to every sequence $\{\ve x^{(k)}\}_{k\in\N}$ and $\{\ve d^{(k)}\}_{k\in\N}$ satisfying the assumptions of Proposition \ref{ProArm}, not only for sequences defined as $\ve x^{(k+1)} = \ve x^{(k)}+\lambda^{(k)}\ve d^{(k)}$. \\
\noindent A further direct consequence of the Armijo condition is the following lemma which will be used in the next section to prove the convergence of the scaled gradient projection method.
\begin{lemma}\label{Lem:1}
Let $\{\ve x^{(k)}\}_{k\in\N}$ be a sequence of points in $\Omega$ and $\{\ve d^{(k)}\}_{k\in\N}$ be a sequence of descent directions such that condition \eqref{A2} holds. Suppose that there exists $l \in \R$ such that $f(\x) \geq l$ for all $\x \in \Omega$ and that
\begin{equation}\label{suff-decr}
f(\xkk) \leq f(\xk+\lamk\dk) \quad \forall k \in \N.
\end{equation}
Then we have
\begin{equation}\label{Lemma:1}
0\leq-\sum_{k=0}^\infty \lamk\nabla f(\xk)^T\dk <\infty.
\end{equation}
\end{lemma}\noindent
{\it Proof.} Inequality \eqref{Armijo} can be rewritten as
\begin{equation*}
-\beta \lamk \nabla f(\xk)^T\dk \leq f(\xk)-f(\xkk).
\end{equation*}
Summing the previous inequality for $k=0,...,j$ gives
\begin{eqnarray}
-\beta \sum_{k=0}^j\lamk\nabla f(\xk)^T\dk &\leq \sum_{k=0}^j (f(\xk)-f(\xkk))\nonumber\\
 &= f(\ve x^{(0)})-f(\xjj)\nonumber\\
 &\leq f(\ve x^{(0)})-l.\label{nuova4}
\end{eqnarray}
Thus, inequality \eqref{Lemma:1} follows.
\endproof
The previous results hold, in general, for any sequence $\{\dk\}_{k\in\N}$ of descent directions. In particular, if we choose as descent direction the vector defined in \eqref{iter0}--\eqref{y1}, then a further property of $\{\dk\}_{k\in\N}$ holds true, as reported in the following lemma whose proof can be found in \cite[Lemmata 2.2,2.3]{Bonettini2009}.
\begin{lemma}
Let $\xk\in\Omega$ and $\dk$ be defined as in \eqref{iter0}--\eqref{y1}. Then we have
\begin{equation}\label{nuovissima1}
\nabla f(\xk)^T\dk \leq -\frac{\|\dk\|^2_{\Dkk}}{\ak}.
\end{equation}
Moreover, $\dk=\ve{0}$ if and only if $\xk$ is stationary for problem \eqref{minf}.
\end{lemma}
We are now ready to give the more general convergence result based on the above mentioned properties of the descent direction and with the Armijo rule establishing the sufficient decrease of the objective function. Its proof is omitted since it can be easily derived by the analogous results in \cite{Birgin2003,Bonettini2009,Bonettini2011}.
\begin{theorem}\label{teo-suff-decr}
Let $\alpha_{\min},\alpha_{\max},\mu$ be three positive constants such that $0<\alpha_{\min}\leq \alpha_{\max}$ and $\mu \geq 1$. Let $\{\ak\}_{k\in\N}\subset [\alpha_{\min},\alpha_{\max}]$ be a sequence of parameters and $\{D_k\}_{k\in\N}\subset {\mathcal M}_\mu$. Let $\{\ve x^{(k)}\}_{k\in\N}\subset\Omega $ be any sequence satisfying property \eqref{suff-decr}, where $\dk$ is defined in \eqref{iter0}--\eqref{y1} and $\lamk$ is computed with the Armijo linesearch procedure in Algorithm \ref{Algo2}. If $\bar{\ve x}$ is a limit point of $\{\xk\}_{k\in\N}$, then $\bar{\ve x}$ is a stationary point for problem \eqref{minf}.
\end{theorem}
It is worth stressing that the previous result applies also to nonconvex problems and the gradient of the objective function is not required to be Lipschitz continuous. Moreover, the only limitations to the algorithms parameters choice are that $\ak$ and the eigenvalues of $D_k$ have to be bounded above and below away from zero.\\
When $\nabla f$ satisfies some Lipschitz property, the next proposition states that the Armijo steplengths are bounded away from zero. This also means that there exists a finite upper bound for the number of backtracking reductions at any iteration (a similar result can be found in \cite[Theorem 3.2]{Auslender2007}). This result will be useful in the convergence rate analysis of the next section.
\begin{proposition}\label{pro:lipschitz}
Assume that $\nabla f$ satisfies one of the following conditions:
\begin{itemize}
\item[a)] $\nabla f$ is globally Lipschitz on $\Omega$;
\item[b)] $\nabla f$ is locally Lipschitz and $f$ is level bounded on $\Omega$.
\end{itemize}
Let $\{\xk\}_{k\in\N}$ be any sequence satisfying the assumptions of Theorem \ref{teo-suff-decr} and $\{\lamk\}_{k\in\N}$ the related steplengths computed by Algorithm \ref{Algo2}. Then, there exists a positive constant $0<\lambda_{\min}\leq 1$ such that
\begin{equation}\label{lambdastar}
\lamk\geq \lambda_{\min}.
\end{equation}
\end{proposition} \noindent
{\it Proof.} If $\nabla f$ is Lipschitz continuous on $\Omega$ with Lipschitz constant $L$, then from the descent lemma \cite[p.667]{Bertsekas1999} we have
\begin{equation}\label{descent-lemma}
f(\xk+\lambda \dk) \leq f(\xk) +\lambda \nabla f(\xk)^T\dk +\frac L 2 \lambda^2\|\dk\|^2,
\end{equation}
where $\lambda\in[0,1]$.\\
If, instead, $\nabla f$ is only locally Lipschitz, by assumption $f$ is level bounded; since \eqref{suff-decr} implies $\{\xk\}_{k\in\N}\subset \Omega_{f(\ve x^{(0)})}$, we have that $\{\xk\}_{k\in\N}$ is bounded. Equation \eqref{nuovissima1} implies $\|\yk-\xk\|\leq \frac 2 {\mu\alpha_{\max}} \|\nabla f(\xk)\|$. Then, $\{\yk\}_{k\in\N}$ is also bounded and there exists a compact set $K$ containing the points $\xk+\lambda(\yk-\xk)$ for any $\lambda\in[0,1]$ and any $k\in\N$. As a consequence of this, inequality \eqref{descent-lemma} holds with $L=L_K$.\\
By inequalities \eqref{descent-lemma} and \eqref{nuovissima1} we further obtain
\begin{eqnarray*}
f(\xk+\lambda \dk) &\leq& f(\xk) +\lambda\nabla f(\xk)^T\dk-\frac L \gamma \lambda^2\nabla f(\xk)^T\dk\\
&=& f(\xk) +\lambda\left(1-\frac L \gamma\lambda\right)\nabla f(\xk)^T\dk,
\end{eqnarray*}
where $\gamma = {\mu\alpha_{\max}}$.
The previous inequality ensures that the Armijo condition
\begin{equation}\label{armijo2}
f(\xk+\lambda\dk)\leq f(\xk)+\lambda\beta\nabla f(\xk)^T\dk\end{equation}
is satisfied, for all $k\in\N$, when $(1-L\lambda/ \gamma)\geq \beta$, that is for all $\lambda$ such that $\lambda \leq (1-\beta)\gamma/L$. If $\lamk$ is the steplength computed by Algorithm \ref{Algo2} and the backtracking loop is performed at least once, then $\lambda =\lamk/\delta$ does not satisfies inequality \eqref{armijo2}, which means $\lamk > \gamma(1-\beta)\delta/L$. Thus, the steplength sequence $\{\lamk\}_{k\in\N}$ satisfies inequality \eqref{lambdastar} with $\lambda_{\min} = \min\{1,\gamma(1-\beta)\delta/L\}$.
\endproof
\section{Convergence analysis of the scaled gradient projection algorithm}\label{SGPconv}
In this section we consider the SGP method whose basic scheme is reported in Algorithm \ref{algo:SGP}. Clearly, Theorem \ref{teo-suff-decr} applies also to Algorithm \ref{algo:SGP}, establishing that any limit point of the sequence $\{\xk\}_{k\in\N}$ is stationary.
Our aim is to propose practical conditions for selecting the SGP metric, i.e. the parameter $\mu_k$, ensuring the convergence of the sequence $\{\xk\}_{k\in\N}$ to a solution of \eqref{minf}, under the only assumption that $f$ is convex and admits a finite minimum.\\
The same conditions allow us also to prove a ${\mathcal O}(1/k)$ complexity result for SGP, which holds when the gradient of $f$ satisfies some Lipschitz assumption.\\

\begin{algorithm}[ht]\caption{Scaled gradient projection (SGP) method}\label{algo:SGP}
Choose $0<\alpha_{\min}\leq\alpha_{\max}$, $\mu \geq 1$, $\delta, \beta \in(0,1)$, $\x^{(0)}\in\Omega$. \\
For $k=0,1,2,...$
\begin{itemize}
\item[1.] Choose $\ak\in[\alpha_{\min},\alpha_{\max}]$;
\item [2.] Choose $\mu_k\leq \mu$ and a positive definite matrix $D_k\in {\mathcal M}_{\mu_k}$;
\item[3.] Compute $\yk$ as in \eqref{y1};
\item[4.] Set $\dk = \yk-\xk$;
\item[5.] Compute the steplength parameter $\lamk$ with Algorithm \ref{Algo2};
\item[6.] Set $\xkk = \xk+\lamk\dk$.
\end{itemize}
\end{algorithm}

\noindent Before to give the main convergence result, we prove the following lemma.
\begin{lemma}\label{lemma:3}
Let $\{\mu_k\}_{k\in\N}$, $\{\zeta_k\}_{k\in\N}$ be two sequences of numbers such that
\begin{equation}\label{condLk}
\mu_k^2 = 1+\zeta_k, \ \ \ \zeta_k \geq 0, \ \ \ \sum_{k=0}^\infty \zeta_k <\infty.
\end{equation}
Then the sequence $\{\theta_k\}_{k \in \N}$, with $\theta_k=\prod_{j=0}^k\mu_j^2$, is bounded.
\end{lemma} \noindent
{\it Proof.} We want to show that there exists a constant $M>0$ such that $\theta_k\leq M$ for all $k\in\N$. By the monotonicity of the logarithm, this is true if and only if $\log(\theta_k)\leq \log(M)$ $\forall k \in \N$. By definition of $\theta_k$ we have
\begin{equation}\label{serie}
\log(\theta_k) = \sum_{j=0}^k \log(\mu_j^2)\leq\sum_{j=0}^\infty \log(\mu_j^2).
\end{equation}
Thus, if the series on the right hand side of \eqref{serie} converges, the quantities $\theta_k$ are bounded for all $k$. We observe that, since $\mu_j^2 = 1+\zeta_j$, by the known limit $\lim_{\zeta_j \rightarrow 0} \log(1+\zeta_j)/\zeta_j = 1$, the series $\sum_{j=0}^\infty \log(\mu_j^2)$ and $\sum_{j=0}^\infty \zeta_j$ have the same behaviour. Thus, since by hypothesis the latter one is convergent, the theorem follows.
\endproof
The next theorem states that, when $f$ is convex and admits finite minimum, if the scaling matrices $D_k$ asymptotically reduce to the identity matrix at a certain rate, then the sequence generated by SGP converges to a solution of \eqref{minf}. The line of the proof is similar to that of \cite[Theorem 1]{Iusem2003}, which can be considered as a special case of it. After giving the proof of our result, we discuss the relations of our approach with the related work already present in the literature.
\begin{theorem}\label{thm:sgp_convex_converge}
Assume that the objective function of \eqref{minf} is convex and the solution set $X^*$ is not empty. Let $\{\xk\}_{k \in \N}$ be the sequence generated by SGP where $D_k\in \mathcal{M}_{\mu_k}$ and $\{\mu_k\}_{k \in \N}$ satisfies \eqref{condLk}. Then the sequence $\{\xk\}_{k \in \N}$ converges to a solution of \eqref{minf}.
\end{theorem} \noindent
{\it Proof.} We recall first the basic norm equality
\begin{equation}\label{norm_equality}
\|\x-\y\|^2_E + \|\y-\z\|^2_E-\|\x-\z\|^2_E = 2(\y-\x)^TE(\y-\z)
\end{equation}
which holds true for any positive definite matrix $E$. Moreover, it is easy to see that if $D_k\in \mathcal{M}_{\mu_k}$, then $D_k^{-1}\in \mathcal{M}_{\mu_k} $. \\
Let $\hatx \in X^*$. By definition of $\yk$ we have
\begin{equation*}
(\yk-\xk+\ak D_k\nabla f(\xk))^TD_k^{-1}(\x-\yk)\geq 0 \ \ \ \forall \x\in\Omega
\end{equation*}
which, for $\x = \hatx$ gives
\begin{eqnarray*}
\fl (\yk-\xk)^T D_k^{-1}(\hatx-\xk) \geq \ak\nabla f(\xk)^T(\xk-\hatx)\\
\qquad \qquad \qquad \quad +(\yk-\xk+\ak D_k\nabla f(\xk))^T\Dkk(\yk-\xk)\\
\geq \ \ak(f(\xk)-f(\hatx))+ \|\yk-\xk\|^2_{\Dkk}+\ak\nabla f(\xk)^T(\yk-\xk)\\
 = \ \ak(f(\xk)-f(\hatx))+\frac{1}{(\lamk)^2}\|\xkk-\xk\|_{\Dkk}^2 + \ak\nabla f(\xk)^T(\yk-\xk),
\end{eqnarray*}
where the inequality follows from the convexity of $f$ and the last equality by definition of $\xkk$. By equality \eqref{norm_equality} with $\x =\xkk$, $\y = \xk$, $\z = \hatx$, $E=\Dkk$ we obtain
\begin{eqnarray*}
\fl \|\xkk-\hatx\|^2_\Dkk = \|\xk-\hatx\|^2_\Dkk+\|\xkk-\xk\|^2_\Dkk-2(\xk-\xkk)^T\Dkk(\xk-\hatx)\\
=\|\xk-\hatx\|^2_\Dkk+\|\xkk-\xk\|^2_\Dkk -2\lamk(\yk-\xk)^T\Dkk(\hatx-\xk)\\
\leq \|\xk-\hatx\|^2_\Dkk+\left(1-\frac 2 {\lamk}\right)\|\xkk-\xk\|^2_\Dkk-2\ak\lamk\nabla f(\xk)^T(\yk-\xk)\\ 
\ \ \ -2\lamk\ak(f(\xk)-f(\hatx))
\end{eqnarray*}
which, since $\lamk \leq 1$, results in
\begin{eqnarray}\label{nuova1}
\|\xkk-\hatx\|^2_\Dkk&\leq &\|\xk-\hatx\|^2_\Dkk-2\ak\lamk\nabla f(\xk)^T(\yk-\xk)+\nonumber\\
&& -2\lamk\ak(f(\xk)-f(\hatx)) \\
&\leq & \|\xk-\hatx\|^2_\Dkk-2\ak\lamk\nabla f(\xk)^T(\yk-\xk). \nonumber
\end{eqnarray}
(since $f(\xk)-f(\hatx)\geq 0$). From the last inequality and in view of \eqref{ine_norm}, it follows that
\begin{eqnarray*}
\frac 1{\mu_k}\|\xkk-\hatx\|^2 &\leq& \|\xkk-\hatx\|^2_\Dkk\\
&\leq&\|\xk-\hatx\|^2_\Dkk-2\ak\lamk\nabla f(\xk)^T(\yk-\xk)\\
&\leq&\mu_k\|\xk-\hatx\|^2-2\ak\lamk\nabla f(\xk)^T(\yk-\xk),
\end{eqnarray*}
that is
\begin{equation*}
\|\xkk-\hatx\|^2\leq \mu_k^2\|\xk-\hatx\|^2-2\mu_k\ak\lamk\nabla f(\xk)^T(\yk-\xk).
\end{equation*}
Recalling that the scalar product at the right-hand-side is nonpositive, since $\mu_k\geq 1$ and $\alpha_k\leq\alpha_{\max}$ this results in
\begin{equation*}
\|\xkk-\hatx\|^2\leq \mu_k^2\|\xk-\hatx\|^2-2{\alpha_{\max}}{\mu_k^2}\lamk\nabla f(\xk)^T(\yk-\xk).
\end{equation*}
By repeatedly applying the previous inequality we obtain
\begin{equation*}
\|\xkk-\hatx\|^2\leq \theta^k_0\|\x^{(0)}-\hatx\|^2-2{\alpha_{\max}}\sum_{j=0}^k{\theta^k_j}\lambda^{(j)}\nabla f(\x^{(j)})^T(\y^{(j)}-\x^{(j)}),
\end{equation*}
where $\theta^k_j = \prod_{i=j}^k\mu_j^2$. Since $\mu_j^2\geq 1$, we have $\theta^k_j\leq \theta_0^k$, and by Lemma \ref{lemma:3} we obtain
\begin{equation}\label{newine}
\|\xkk-\hatx\|^2\leq M\|\x^{(0)}-\hatx\|^2-2{\alpha_{\max}}M\sum_{j=0}^k\lambda^{(j)}\nabla f(\x^{(j)})^T(\y^{(j)}-\x^{(j)}),
\end{equation}
where $\theta_0^k\leq M$. Now we can apply Lemma \ref{Lem:1} to conclude that $\{\xk\}_{k \in \N}$ is bounded and, thus, it has at least one limit point. Let us denote such limit point by $\x^{\infty}$. By Theorem \ref{teo-suff-decr}, $\x^{\infty}$ is stationary; in particular, since $f$ is convex, it is a minimum point, i.e. $\x^{\infty}\in X^*$. Let $\{\x^{(k_i)}\}_{i \in \N}$ be a subsequence of $\{\xk\}_{k \in \N}$ which converges to $\x^\infty$. By applying the same arguments employed to derive \eqref{newine}, for any fixed $i \in \N$ and for all $k\geq k_i$ we obtain
\begin{equation}\label{boh}
\|\xk-\x^\infty\|^2\leq M\|\x^{(k_i)}-\x^\infty\|^2-2{\alpha_{\max}}M\sum_{j=k_i}^k\lambda^{(j)}\nabla f(\x^{(j)})^T(\y^{(j)}-\x^{(j)}).
\end{equation}
Since $\{\x^{(k_i)}\}_{i \in \N}$ converges to $\x^\infty$ and $-\sum_{j=0}^\infty\lambda^{(j)}\nabla f(\x^{(j)})^T(\y^{(j)}-\x^{(j)})$ is a convergent series, for any $\varepsilon > 0$ there exists a sufficiently large integer $k_i$ such that $\|\x^{(k_i)}-\x^\infty\|^2\leq \varepsilon/2M$ and $-\sum_{j=k_i}^k\lambda^{(j)}\nabla f(\x^{(j)})^T(\y^{(j)}-\x^{(j)})\leq \varepsilon/(4M\alpha_{\max})$. Then, it follows from \eqref{boh} that $\|\xk-\x^\infty\|^2\leq \varepsilon$ for all $k\geq k_i$. Since $\varepsilon$ can be chosen arbitrarily small, this means that the whole sequence $\{\xk\}_{k \in \N}$ converges to $\x^{\infty}$.
\endproof
The previous theorem gives an easily implementable rule to ensure the theoretical convergence of SGP to a solution. Moreover, as shown in section \ref{sec:num}, it seems to have a favourable impact also on the practical performances of the method. This result is also coherent with the conclusions drawn from the numerical experience in \cite{Bonettini2013}, where the advantages of using a scaling matrix multiplying the gradient were observed mainly at the initial iterations.\\
Finally, we observe that methods employing a variable scaling are analyzed also in two very recent papers \cite{Combettes2014,Combettes2013} in the context of more general variational problems. In these papers, the authors also analyze the convergence of a variable metric forward--backward algorithm which applies to the convex optimization problem
\begin{equation}\label{minfg}
\min_{\x\in\R^n} f(\x) + g(\x)
\end{equation}
and can be described by the following iteration
\begin{equation}\label{SFB}
\xkk = \xk + \lambda_k({\rm{prox}}_{\alpha_k g}^{D_k}(\xk-\alpha_k D_k\nabla f(\xk)) - \xk),
\end{equation}
where
\begin{equation}\nonumber
{\rm{prox}}_{\alpha_k g}^{D_k}(\y) = { \argmin{\x\in\R^n}}\ g(\x) + \frac 1 {2\alpha_k} (\x-\y)^TD_k^{-1}(\x-\y).
\end{equation}
Clearly, when $g$ is the indicator function of the convex set $\Omega$, problem \eqref{minf} is equivalent to \eqref{minfg} and the SGP iteration can be expressed in the same form of \eqref{SFB}. In \cite{Combettes2014}, the convergence of the iterates \eqref{SFB} is proved for objective functions with Lipschitz continuous gradients, under the condition
\begin{equation}\label{Combettes-cond}
(1+\zeta_k)D_{k+1} \succeq D_k
\end{equation}
where $\zeta_k$ is a summable sequence.\\
Variable metrics were considered also in \cite[Chapter 5]{Nedic2002} in the context of subgradient methods for nonsmooth, convex, unconstrained minimization. In this case, setting $D_k = B_kB_k^T$, the scaling matrices are assumed to satisfy $\prod_{k=0}^\infty\|B_{k+1}^{-1}B_k\|^2<\infty$ and
\begin{equation}\label{Nedic-cond}
\|B_{k+1}^{-1}B_k\|\geq 1.
\end{equation}
We remark that our condition, $D_k\in {\mathcal M}_{\mu_k}$, is quite different from both \eqref{Combettes-cond} and \eqref{Nedic-cond} since it does not impose a strict connection between the scaling matrices at two successive iterates. This freedom of choosing the metric at each iteration allows for example to adopt a suitable adaptation of a well performing scaling technique, based on a gradient splitting \cite{Bertero2008,Lanteri2002}, which may lead to significant improvements of the convergence behaviour, as we will show in section \ref{sec:num}.\\
In the following we give a complexity result about SGP, showing that it has a ${\mathcal O}(1/k)$ convergence rate on the objective function value. Similar results can be found in \cite{Beck2009b} for forward--backward methods with linesearch along the projection arc (i.e. of the form \eqref{SFB} with $D_k=I$, $\lambda_k=1$ for all $k$ and with $\alpha_k$ determined by a backtracking procedure).
\begin{theorem}\label{thm:SGP_convergence_rate}
Assume that the hypotheses of Theorem \ref{thm:sgp_convex_converge} hold and, in addition, that assumption a) or b) of Proposition \ref{pro:lipschitz} is satisfied. Let $f^*$ be the optimal function value for problem \eqref{minf}. Then, we have
\begin{equation*}
f(\xk)-f^* = {\mathcal O}(1/k).
\end{equation*}
\end{theorem}
{\it Proof.} Setting $a = 2\lambda_{\min}\alpha_{\min}$, where $\lambda_{\min}$ is defined in Proposition \ref{pro:lipschitz}, from \eqref{nuova1} we have
\begin{eqnarray*}
\|\xkk-\hatx\|^2_\Dkk &\leq &\|\xk-\hatx\|^2_\Dkk-2\ak\lamk\nabla f(\xk)^T(\yk-\xk) +\\
& &  -2\lamk\ak(f(\xk)-f(\hatx))\\
&\leq& \|\xk-\hatx\|^2_\Dkk -2\alpha_{\max}\lamk\nabla f(\xk)^T(\yk-\xk)  +\\
& & +a(f(\hatx)-f(\xk)),
\end{eqnarray*}
where the second inequality follows from the fact that $\nabla f(\xk)^T(\yk-\xk)$ and $f(\hatx)-f(\xk)$ are negative quantities.
Thanks to inequality \eqref{ine_norm}, we can write
\begin{eqnarray*}
\frac{1}{\mu_k}\|\xkk-\hatx\|^2&\leq&\|\xkk-\hatx\|^2_\Dkk\\
&\leq &\|\xk-\hatx\|^2_\Dkk -2\alpha_{\max}\lamk\nabla f(\xk)^T(\yk-\xk)+\\ & &  +a(f(\hatx)-f(\xk))\\
&\leq &\mu_k\|\xk-\hatx\|^2 -2\alpha_{\max}\lamk\nabla f(\xk)^T(\yk-\xk)+\\ & &  +a(f(\hatx)-f(\xk)).
\end{eqnarray*}
By multiplying the last inequality by $\mu_k$ we obtain
\begin{eqnarray*}
\|\xkk-\hatx\|^2&\leq& \mu_k^2\|\xk-\hatx\|^2 -2\alpha_{\max}\mu_k\lamk\nabla f(\xk)^T(\yk-\xk)+\\ & &  +\mu_k a(f(\hatx)-f(\xk))\\
&\leq& \mu_k^2\|\xk-\hatx\|^2 -2\alpha_{\max}\mu_k^2\lamk\nabla f(\xk)^T(\yk-\xk)+\\ & &  + a(f(\hatx)-f(\xk)),
\end{eqnarray*}
where the last inequality follows from the fact that $\mu_k\geq 1$. By repeatedly applying the last inequality we obtain
\begin{eqnarray}
\|\xkk-\hatx\|^2 &\leq& \theta_0^k \|\x^{(0)}-\hatx\|^2 -2\alpha_{\max}\sum_{j=0}^k\theta^k_j\lamj\nabla f(\xj)^T(\yj-\xj)+ \nonumber\\
& &+ a((k+1) f(\hatx)-\sum_{j=0}^k f(\x^{(j)}))\nonumber\\
&\leq& M \|\x^{(0)}-\hatx\|^2 -2\alpha_{\max}M\sum_{j=0}^k\lamj\nabla f(\xj)^T(\yj-\xj)+ \nonumber\\
& &+ a((k+1) f(\hatx)-\sum_{j=0}^k f(\x^{(j)})),\label{nuova2}
\end{eqnarray}
where, as in the proof of Theorem \ref{thm:sgp_convex_converge}, we set $\theta^k_j = \prod_{i=j}^k\mu_j^2$ and $M$ is the upper bound of all $\theta^k_j$. Thanks to inequality \eqref{nuova4}, we have
\begin{eqnarray}
\|\xkk-\hatx\|^2 &\leq& M \|\x^{(0)}-\hatx\|^2 +\frac{2\alpha_{\max}M}{\beta} (f(\x^{(0)})-f(\hatx))+ \nonumber\\
& &+ a(k f(\hatx)-\sum_{j=1}^k f(\x^{(j)})),\label{nuova3}
\end{eqnarray}
where we also added the positive quantity $a(f(\x^{(0)}) - f(\hatx))$ to the right hand side of \eqref{nuova2}. Moreover, exploiting the inequality
\begin{equation}\nonumber
0\leq \sum_{j=0}^k j (f(\xj)-f(\xjj)) = \sum_{j=1}^k f(\xj) - kf(\xkk)
\end{equation}
gives
\begin{eqnarray*}
\|\xkk-\hatx\|^2 &\leq& M \|\x^{(0)}-\hatx\|^2 +\frac{2\alpha_{\max}M}{\beta} (f(\x^{(0)})-f(\hatx))+ \\
& &+ ak( f(\hatx)- f(\xkk)).
\end{eqnarray*}
Rearranging terms, this finally yields
\begin{equation*}
f(\xkk) - f(\hatx)\leq \frac{M}{ak}\left( \|\x^{(0)}-\hatx\|^2 + 2\frac{\alpha_{\max}}{\beta}(f(\x^{(0)})-f(\hatx))\right),
\end{equation*}
establishing the result.
\endproof
In the recent literature, several authors developed the so-called \emph{intertial} methods, which are first order methods including an extrapolation step which allows to prove a ${\mathcal O}(1/{k^2})$ convergence rate on the objective function values (see for example \cite{Beck2009b,Nesterov2005,Villa2013,Ochs2014}). However, as we will show in section \ref{sec:num}, the practical performances of SGP can be comparable with those of ${\mathcal O}(1/{k^2})$ methods, even if the theoretical convergence rate estimate is only ${\mathcal O}(1/k)$.
\section{Numerical illustration}\label{sec:num}
In this section we consider some relevant applications and we show that they can be effectively solved by algorithms which can be framed in the analysis of the previous sections. We give also some hints on how to choose the parameters $\alpha_k$ and $D_k$ at each iteration, even if a specific treatment of this issue is far beyond the scope of this paper. Both sets of numerical tests concern the image deconvolution problem in the presence of Poisson noise. In particular, in the next subsection we will consider a fit-to-data + regularization model with an arbitrarily fixed regularization parameter, while in the following tests we will investigate the same problem combined with an automatic procedure for the choice of this parameter recently proposed by Zanni et al. \cite{Zanni2015}.
\subsection{Edge preserving image restoration}
\begin{figure}[ht]
\begin{center}
\begin{tabular}{ccc}
\includegraphics[width = 0.3\linewidth]{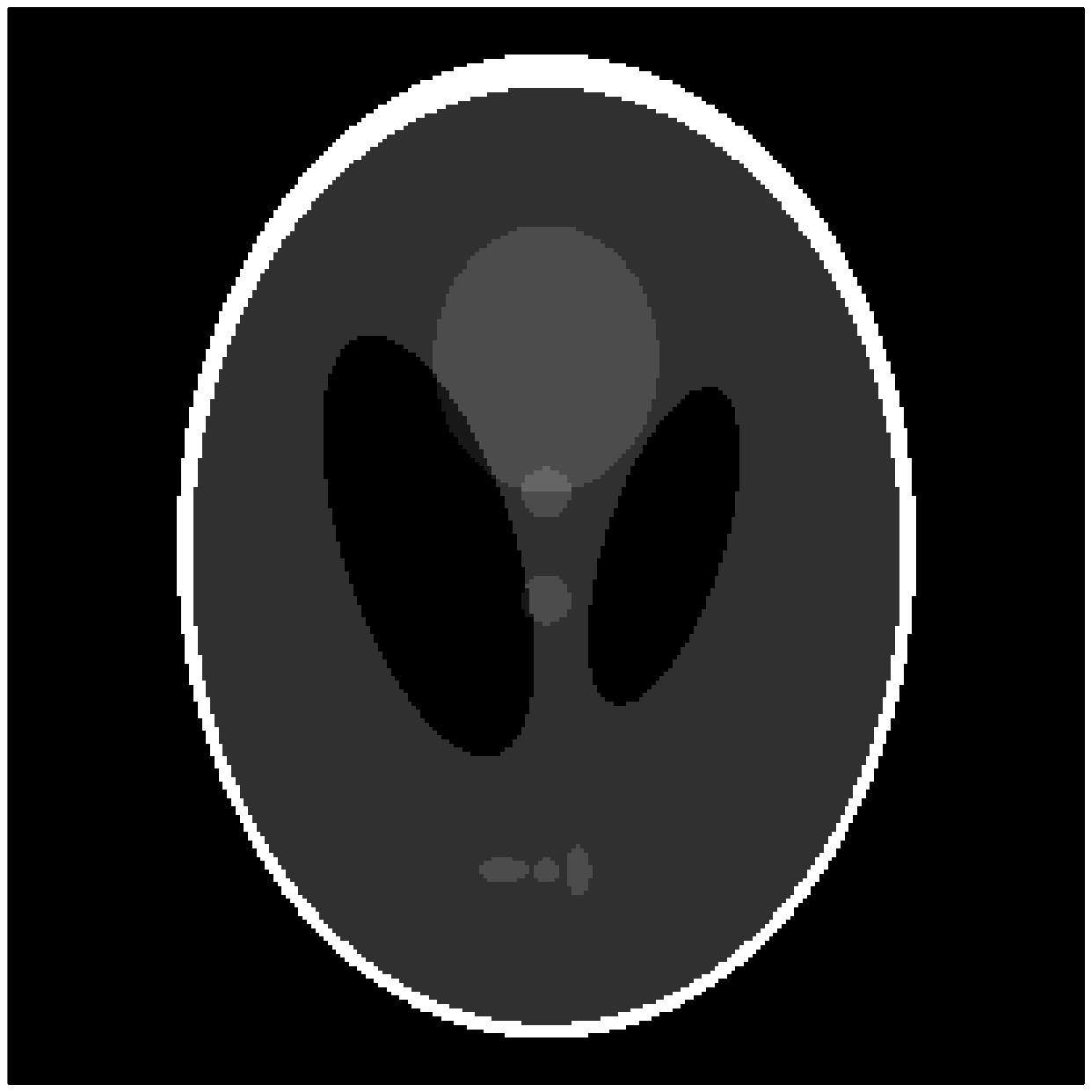} &
\includegraphics[width = 0.3\linewidth]{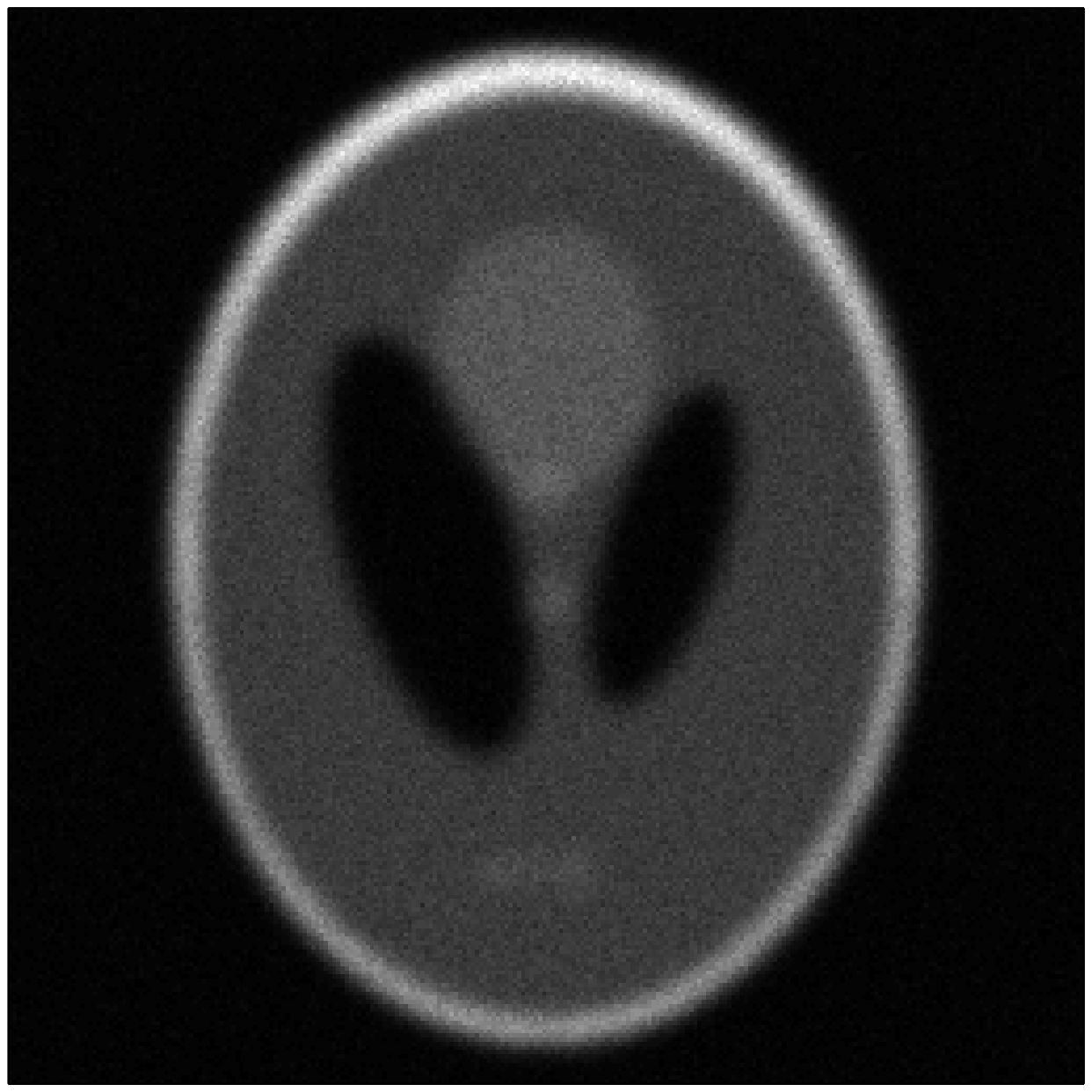} &
\includegraphics[width = 0.3\linewidth]{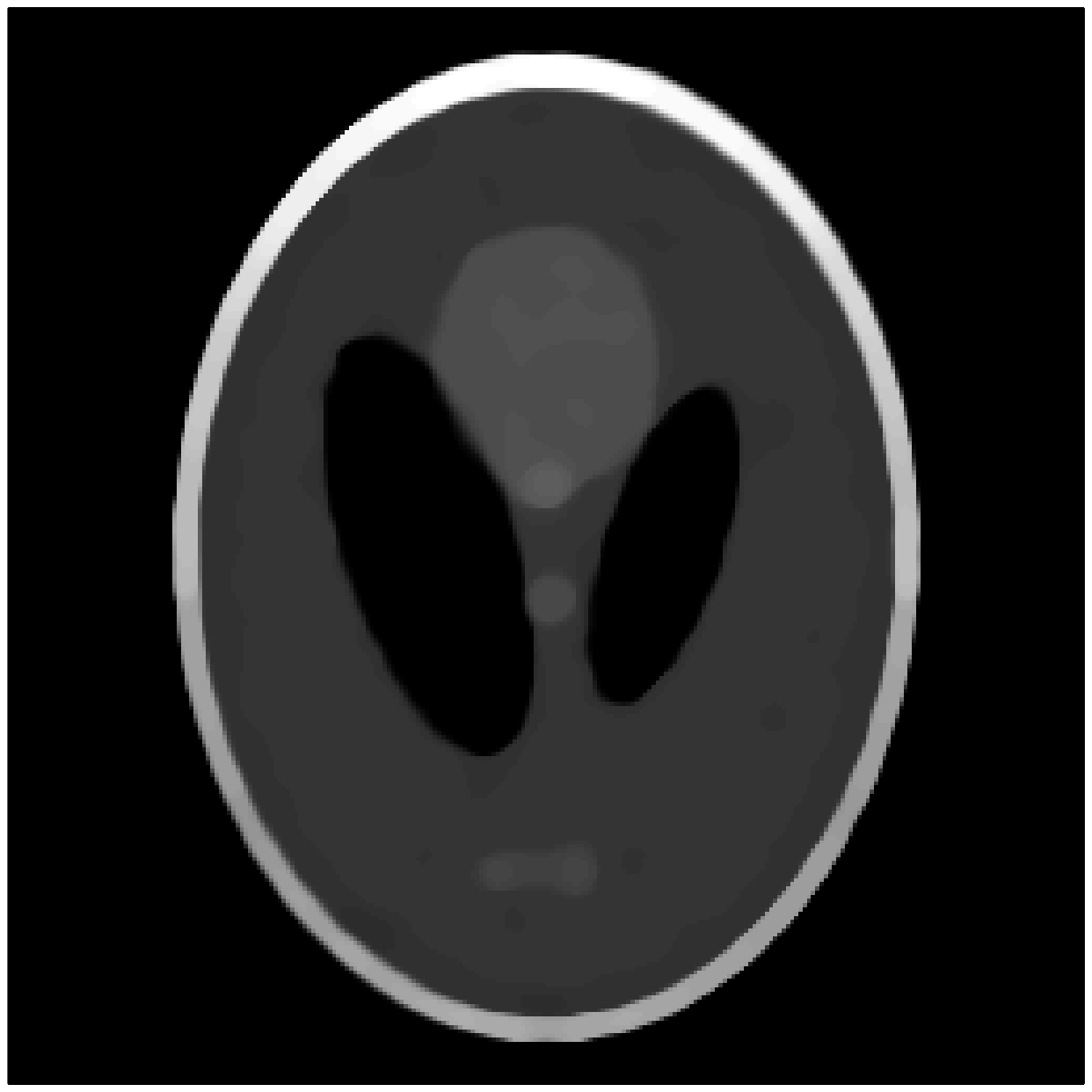}\\
\includegraphics[width = 0.3\linewidth]{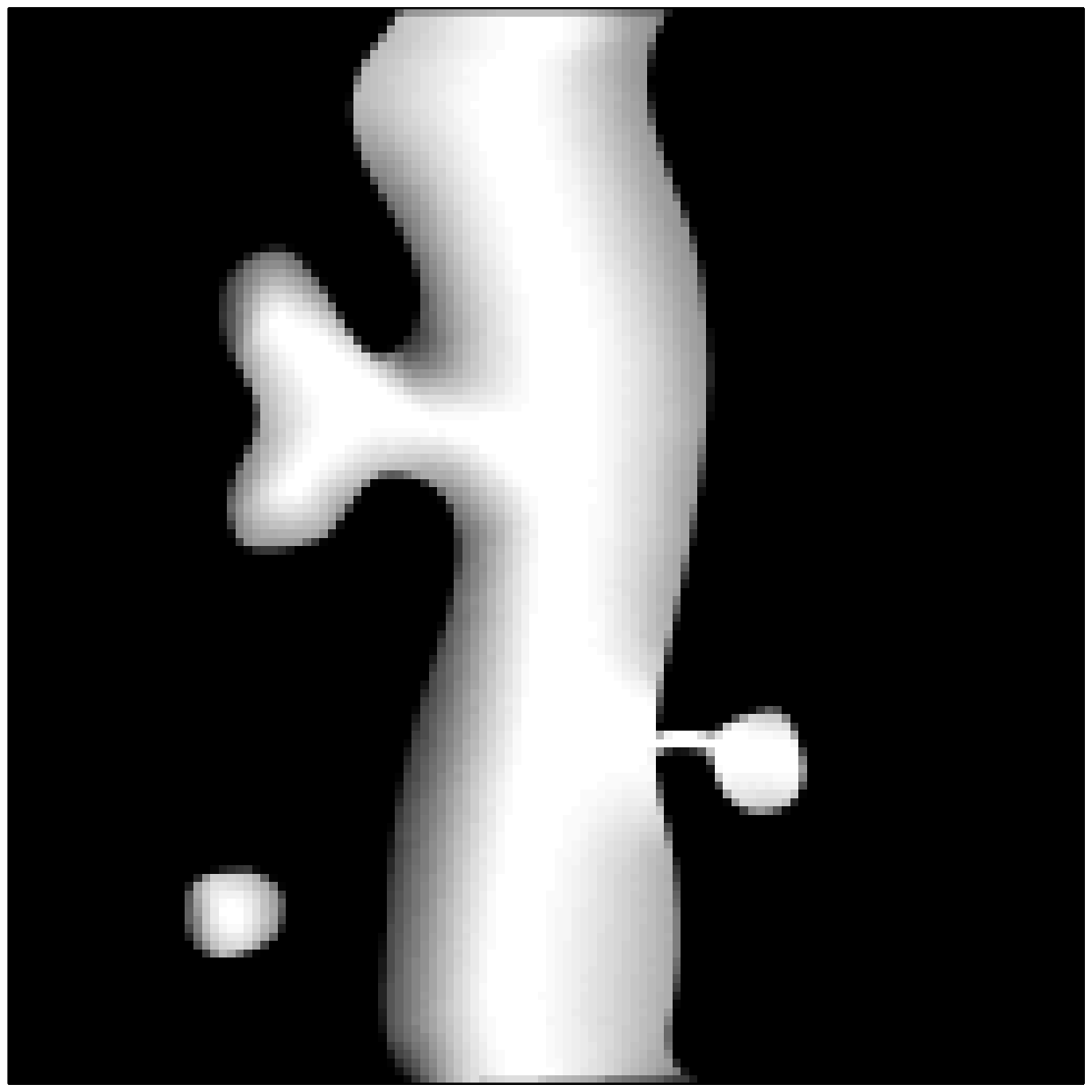} &
\includegraphics[width = 0.3\linewidth]{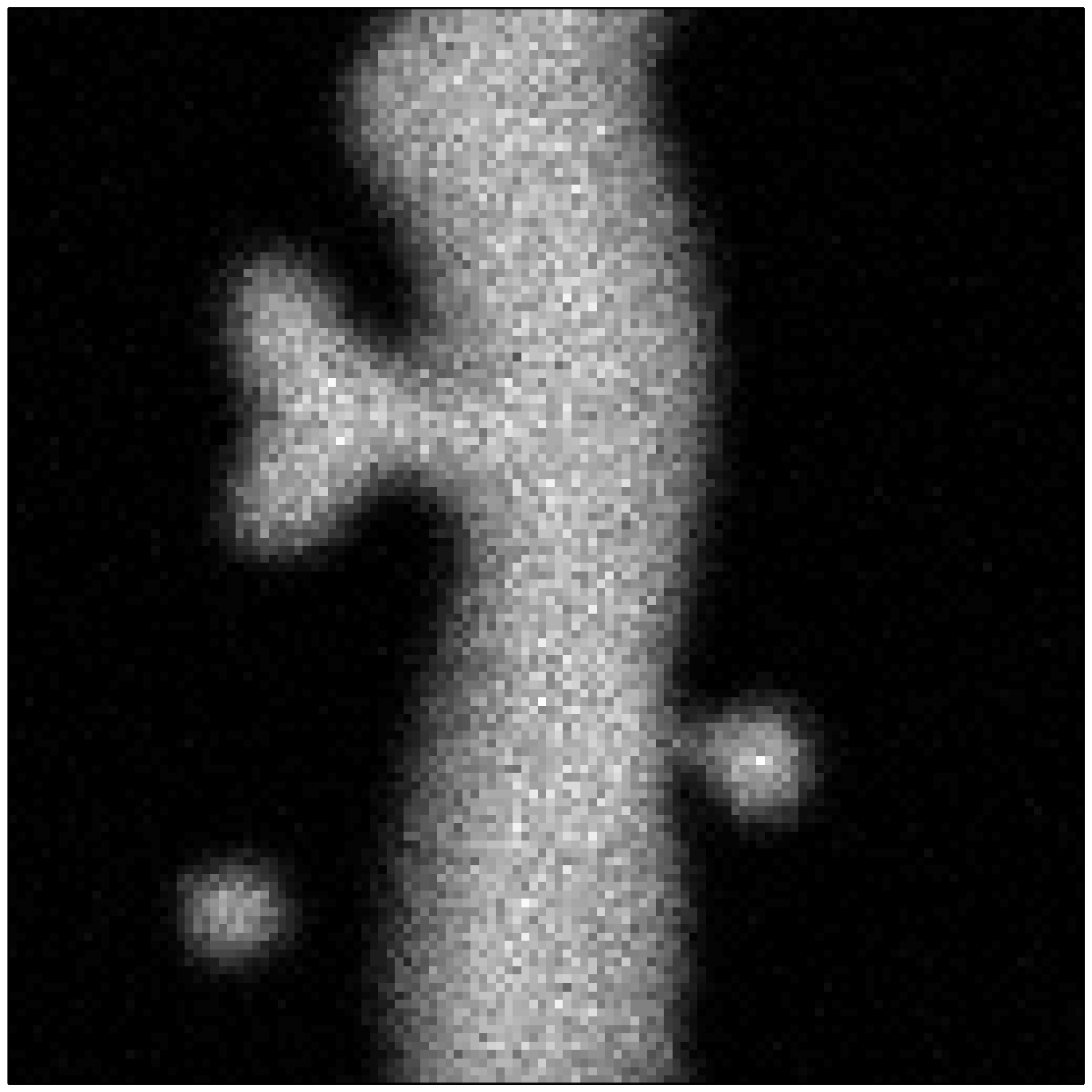} &
\includegraphics[width = 0.3\linewidth]{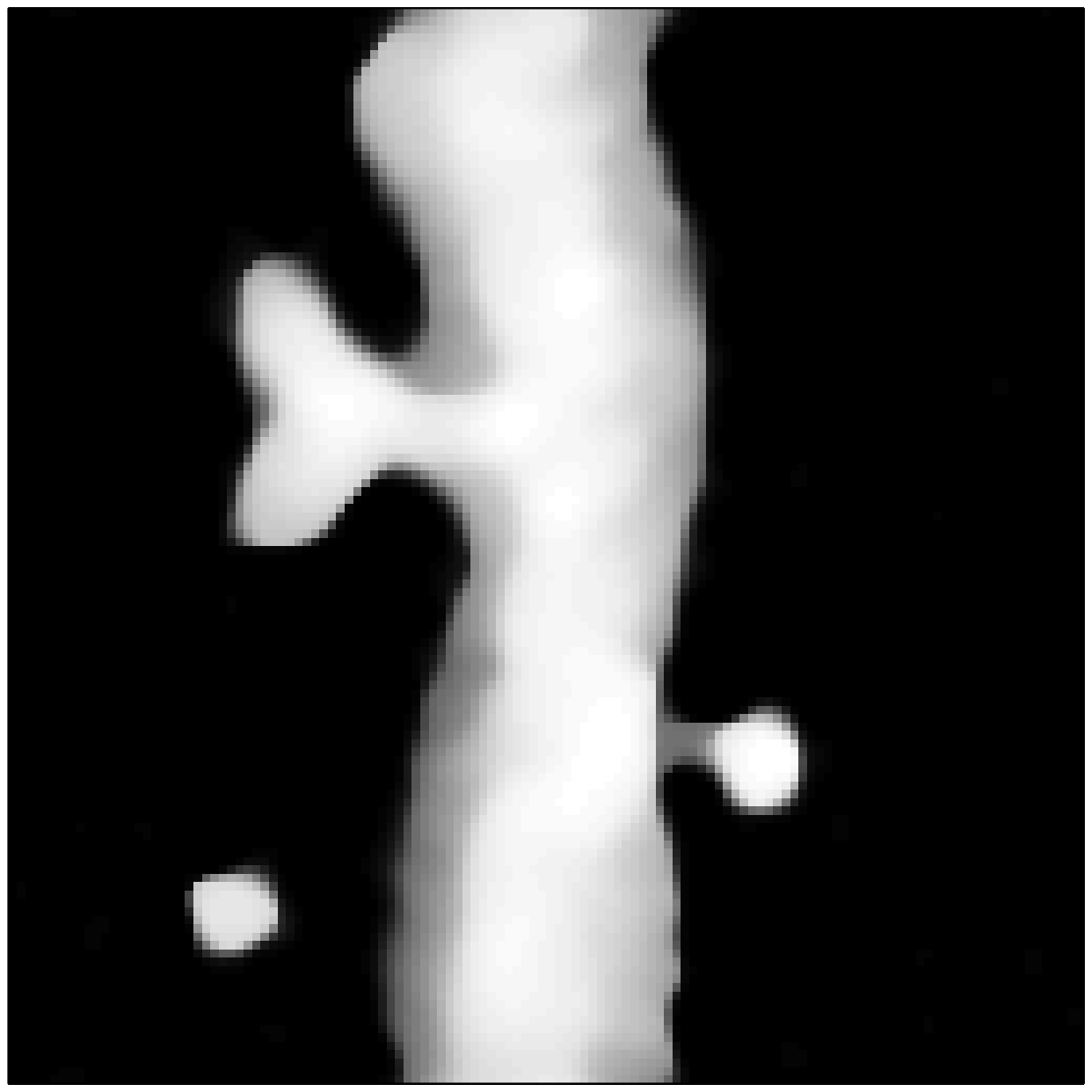}
\end{tabular}
\end{center}
\caption{Test problems: Shepp-Logan (top) and confocal microscopy (bottom) phantoms. Original image (left), noisy blurred image (middle) and optimal solution (right).}
\label{fig:1}
\end{figure}

\begin{figure}[ht]
\begin{center}
\begin{tabular}{cc}
\includegraphics[width = 0.45\linewidth]{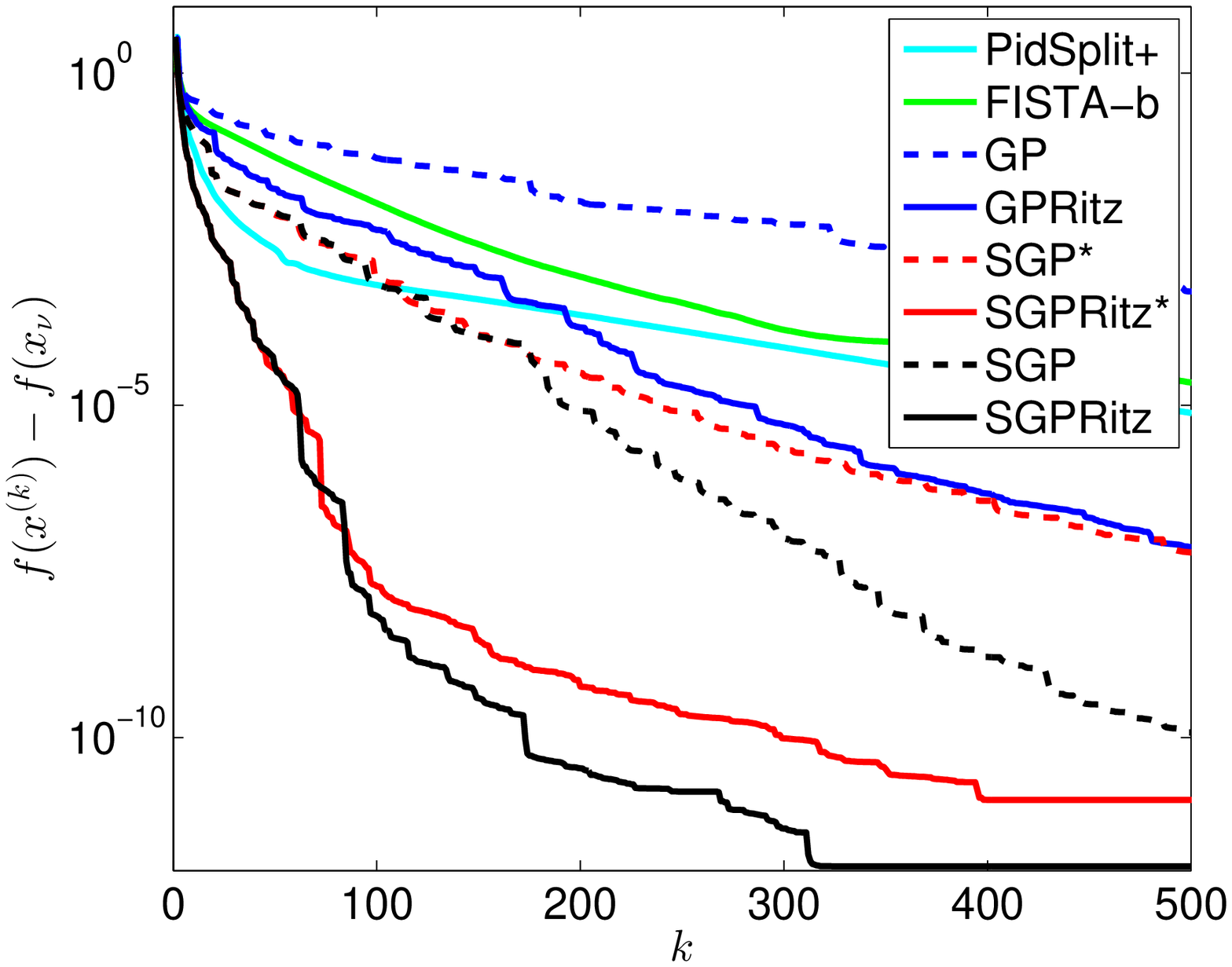} & 
\includegraphics[width = 0.45\linewidth]{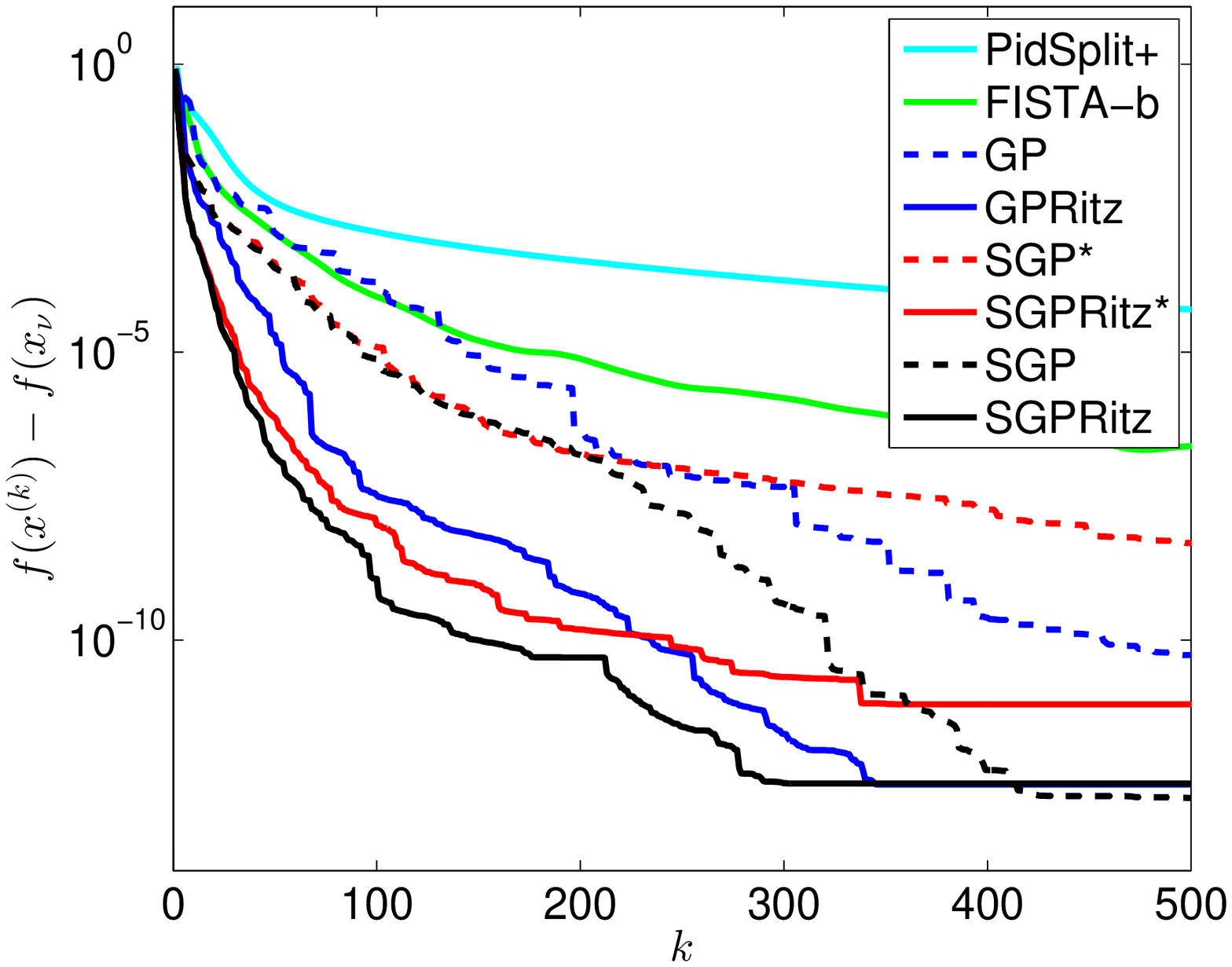}
\end{tabular}
\end{center}
\caption{Image deconvolution: objective function decrease in logarithmic scale versus the iterations number for the SL (left) and CM (right) datasets.}
\label{fig:2}
\end{figure}

Our basic assumption is that the available data $\ve g\in \R^n$ is a realization of a Poisson random variable whose mean is $A\xstar + b\ve e$, where $A\in \R^{n\times n}$ is a structured matrix representing the convolution operator, $b\in \R$ is a positive parameter representing the background radiation, $\ve e\in \R^n$ is the vector of all ones and $\xstar$ is the image we would like to recover. In the following, we will assume that $A\ve e=\ve e$, $A^T\ve e=\ve e$, which is not a restrictive assumption, since it can be assured by a simple normalization.\\
According to the Bayesian approach \cite{Geman1984}, an approximation $\x_\nu$ of $\xstar$ can be obtained by solving the following optimization problem
\begin{equation}\label{KL_R}
\min_{\x\geq 0} f(\x)\equiv KL(\x)+\nu R(\x),
\end{equation}
where $KL(\x)$ is the generalized Kullback--Leibler divergence
\begin{equation}\label{KL}
KL(\x) = \sum_{i=1}^n \left\{g_i\log\left(\frac{g_i}{ (A\x)_i + b}\right)+(A\x)_i + b -g_i \right\},
\end{equation}
$R(\x)$ is some regularization functional, chosen according to the a priori information on the desired solution, and $\nu>0$ is the regularization parameter balancing the relative weight of the two terms. In order to preserve the edges in the restored image, a good choice for the regularization term is the following hypersurface (HS) functional \cite{Acar2004,Bertero2010}
\begin{equation}\label{TV}
HS_\rho(\x) = \sum_{i=1}^n \sqrt{ (\ve D_i^h \x)^2 + (\ve D_i^v\x)^2 + \rho^2},
\end{equation}
where $\rho>0$ and $\ve D_i^h$, $\ve D_i^v$ are finite difference approximations of the horizontal and vertical image gradient, respectively. If $\rho$ is small, it can be considered as an approximation of the total variation functional, but it has been shown that better reconstructions can be obtained for large values of the smoothing parameter \cite{Bonettini2011b}. Thus, we consider the following convex optimization problem
\begin{equation}\label{KLHS}
\min_{\x\geq 0} f(\x)\equiv KL(\x)+\nu HS_\rho(\x),
\end{equation}
whose main features have been studied in \cite{Bonettini2010c}.\\
As for the SGP method, borrowing the ideas in \cite{Zanella2009}, at each iteration $k$ we adopt the following diagonal scaling matrix
\begin{equation}\label{D_k}
[D_k]_{ii} = \max\left\{\frac 1{\mu_k}, \min\left\{\mu_k,\frac{\xk}{1+\nu V^R_i(\xk)}\right\}\right\}, \qquad i=1,\ldots,n,
\end{equation}
where $V^R_i(\xk)$ is defined as in \cite[formula (25)]{Zanella2009}, while $\mu_k = \sqrt{1+10^{10}/k^2}$ so that Theorem \ref{thm:sgp_convex_converge} applies. 
The steplength parameter $\alpha_k$ is then computed in two different ways:
\begin{itemize}
\item the adaptive alternation of the scaled Barzilai--Borwein (BB) rules as proposed in \cite{Bonettini2009};
\item the Ritz-like values proposed by Fletcher \cite{Fletcher2012} for a steepest descent method in the case of unconstrained optimization and recently extended to the SGP algorithm applied to a general constrained problem \eqref{minf} \cite{Porta2015a,Porta2015b}.
\end{itemize}
Besides SGP, we consider also for comparison the ``plain'' gradient projection (GP) method with Euclidean projection and variable steplength (chosen with the same two rules exploited in the scaled case), the PidSplit+ algorithm \cite{Setzer2010}, which is an alternating direction method of multipliers specific for the minimization of the Kullback--Leibler plus the discrete total variation functional ($\rho = 0$), adapted to the smoothed case with $\rho > 0$, and the accelerated proximal-gradient method with inertial/extrapolation with backtracking (FISTA-b) \cite{Beck2009b}.\\
As test problems, we consider:
\begin{itemize}
\item the Shepp-Logan (SL) phantom of size $256\times256$, multiplied by a factor of 500, corrupted with Gaussian blur of variance 9 and with Poisson noise simulated using the \verb"imnoise" Matlab function on the blurred image including the additive background. The background constant is $b=10$;
\item the confocal microscopy (CM) phantom of size $128 \times 128$ described in \cite[section V.C]{Willet2003}, with values in the range $[0,68]$ and with a constant background $b=1$.
\end{itemize}
We assume periodic boundary conditions, so that the matrix $A$ is block circulant with circulant blocks (BCCB) and the matrix-vector products involving $A$ can be performed with a $\mathcal{O}(n\log(n)) $ complexity by means of the fast Fourier transform \cite{Hansen2006}. In figure \ref{fig:1} we report the original objects, the corrupted images and the solutions $\xstar$ of problem \eqref{KL_R} for both test problems. The parameters $(\nu,\rho)$ in \eqref{KLHS} have been empirically tuned to obtain a visually satisfactory solution and have been set equal to $(0.0415,1)$ for SL and $(0.06,1)$ for CM. Moreover, the `$\gamma$' parameter of PidSplit+ has been set equal to $50/\nu$ (SL) and $1/\nu$ (CM) and the initial steplength parameter for FISTA-b is 100 in both cases.
In order to illustrate the convergence behaviour of the methods, we first compute a ground truth solution $\x_\nu$ (see figure \ref{fig:1}, right panel) by running 1500 iterations of SGP. Then, we evaluate the progress towards this solution by computing at each iterate the relative difference of the objective function value with respect to the estimated minimum $f(\x_\nu)$ (see figure \ref{fig:2}). We include in our comparison also the version of SGP with fixed bounds on the scaling matrix $\mu_k=\mu=10^5$, which is denoted by SGP$^*$.
From figure \ref{fig:2} we can observe what follows:
\begin{itemize}
\item the choice of a suitable projection operator can have a significant impact on the practical performances of the gradient projection methods, since GP is outperformed by SGP with both choices for the steplength parameters;
\item SGP with variable bounds on the scaling matrix gives the best performances: in particular, condition \eqref{condLk}, which is employed in Theorem \ref{thm:sgp_convex_converge} to prove the convergence of the method on convex problems, seems also to significantly improve its practical performances, especially when the iterates are close to the solution;
\item in spite of the theoretical convergence rate given in Theorem \ref{thm:SGP_convergence_rate}, the practical behaviour of SGP is comparable with the ${\mathcal O} (1/{k^2})$ method FISTA-b.
\end{itemize}
\subsection{Automatic parameter estimation}
\begin{figure}[h]
\begin{center}
\begin{tabular}{cc}
\includegraphics[width = 0.3\linewidth]{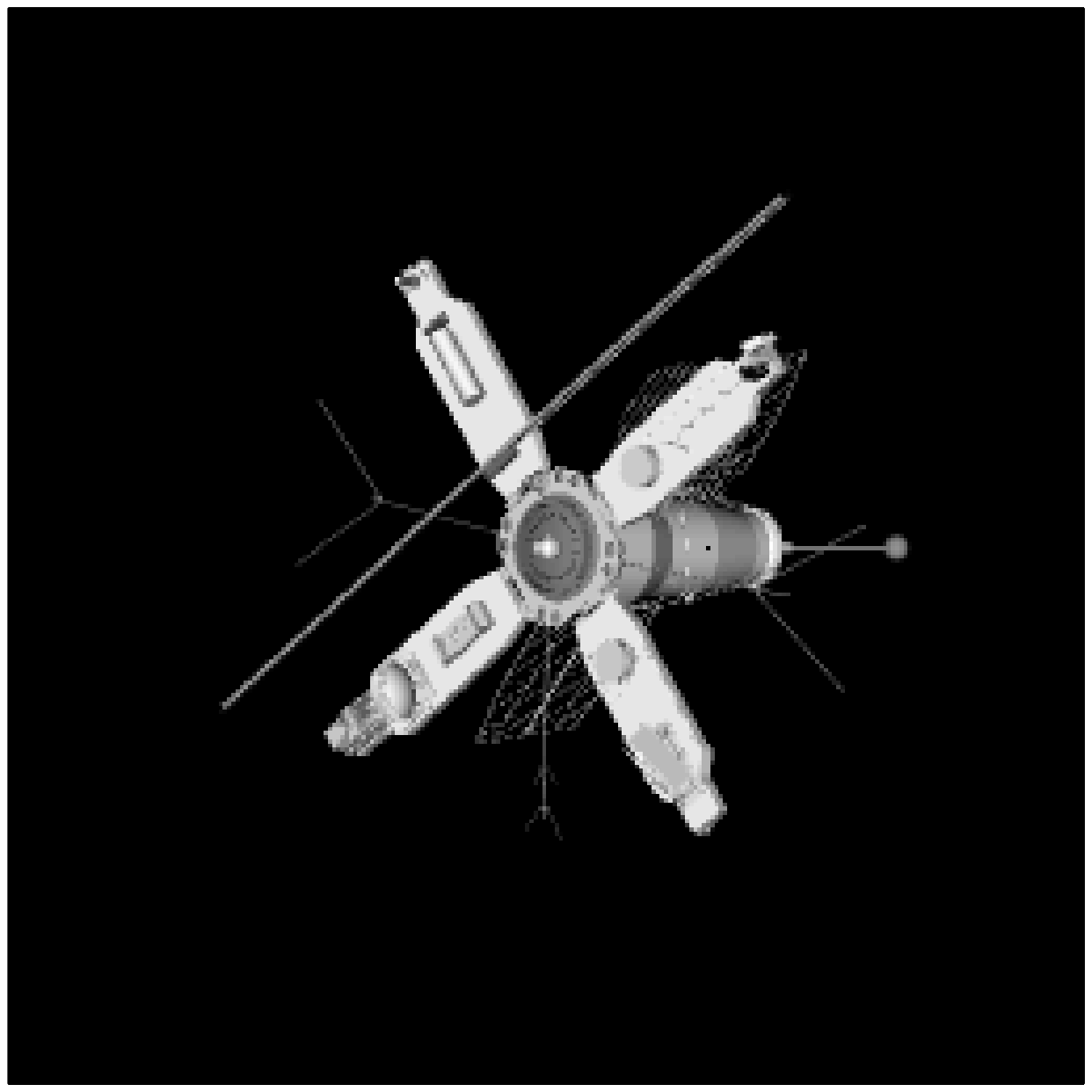} &
\includegraphics[width = 0.3\linewidth]{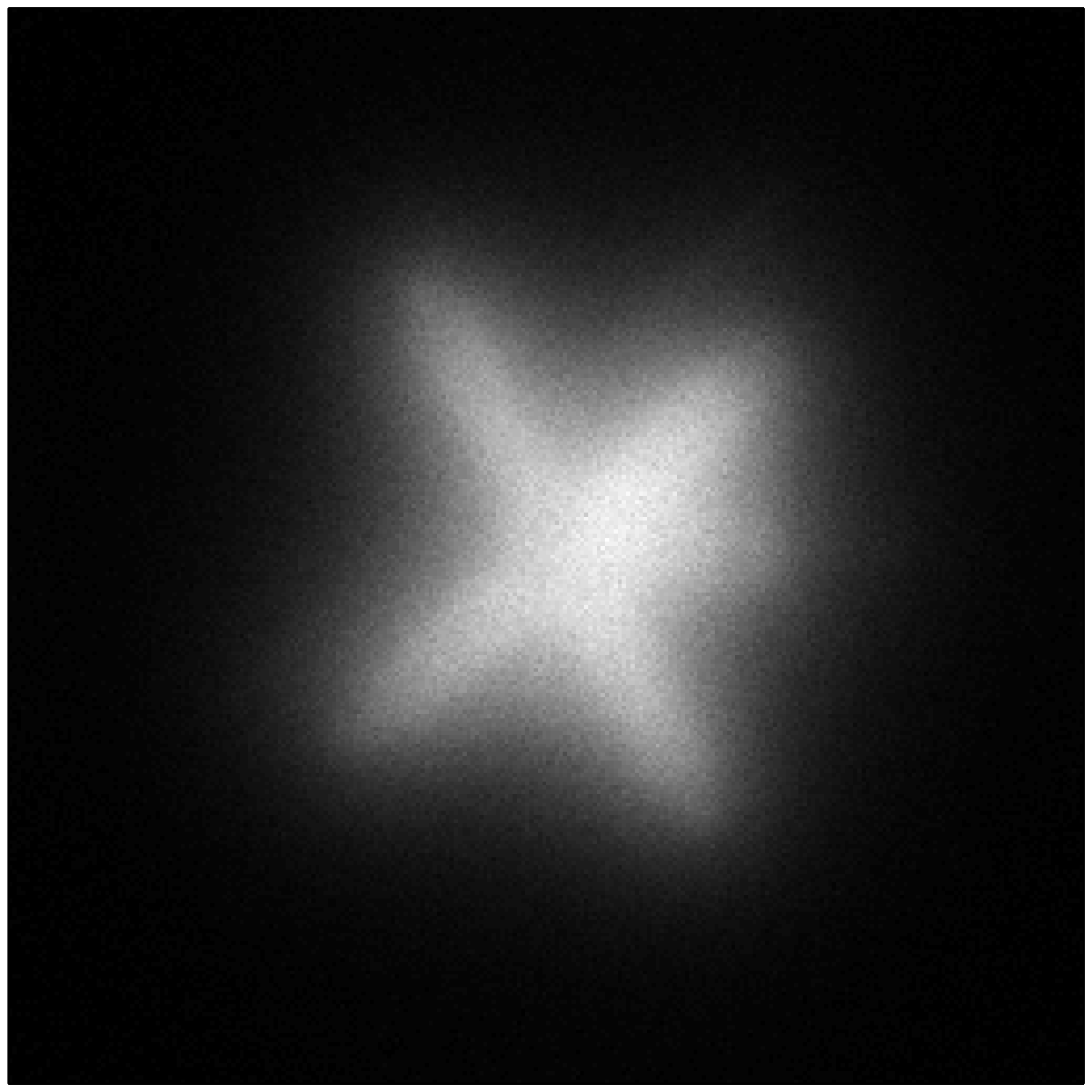}\\
\includegraphics[width = 0.3\linewidth]{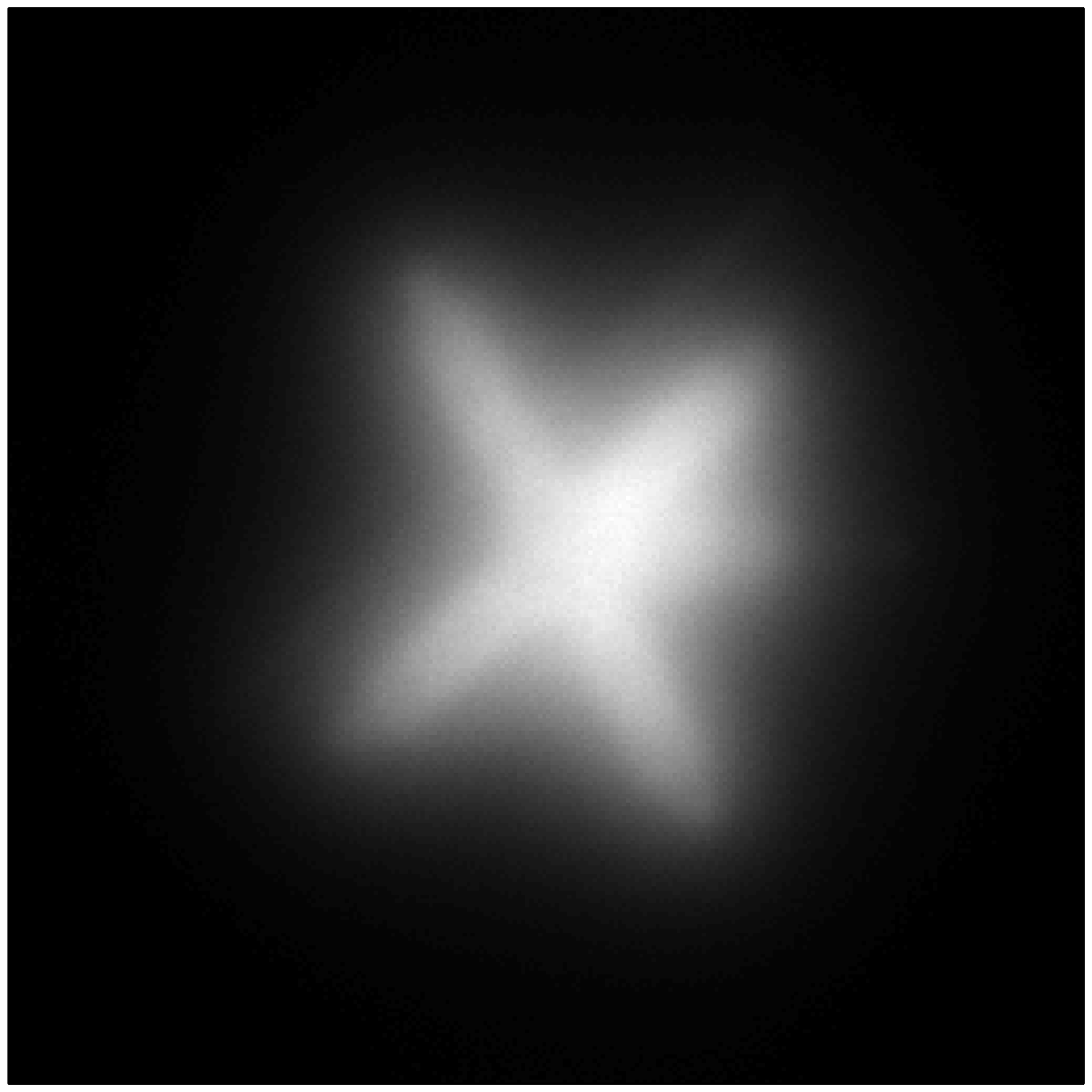} &
\includegraphics[width = 0.3\linewidth]{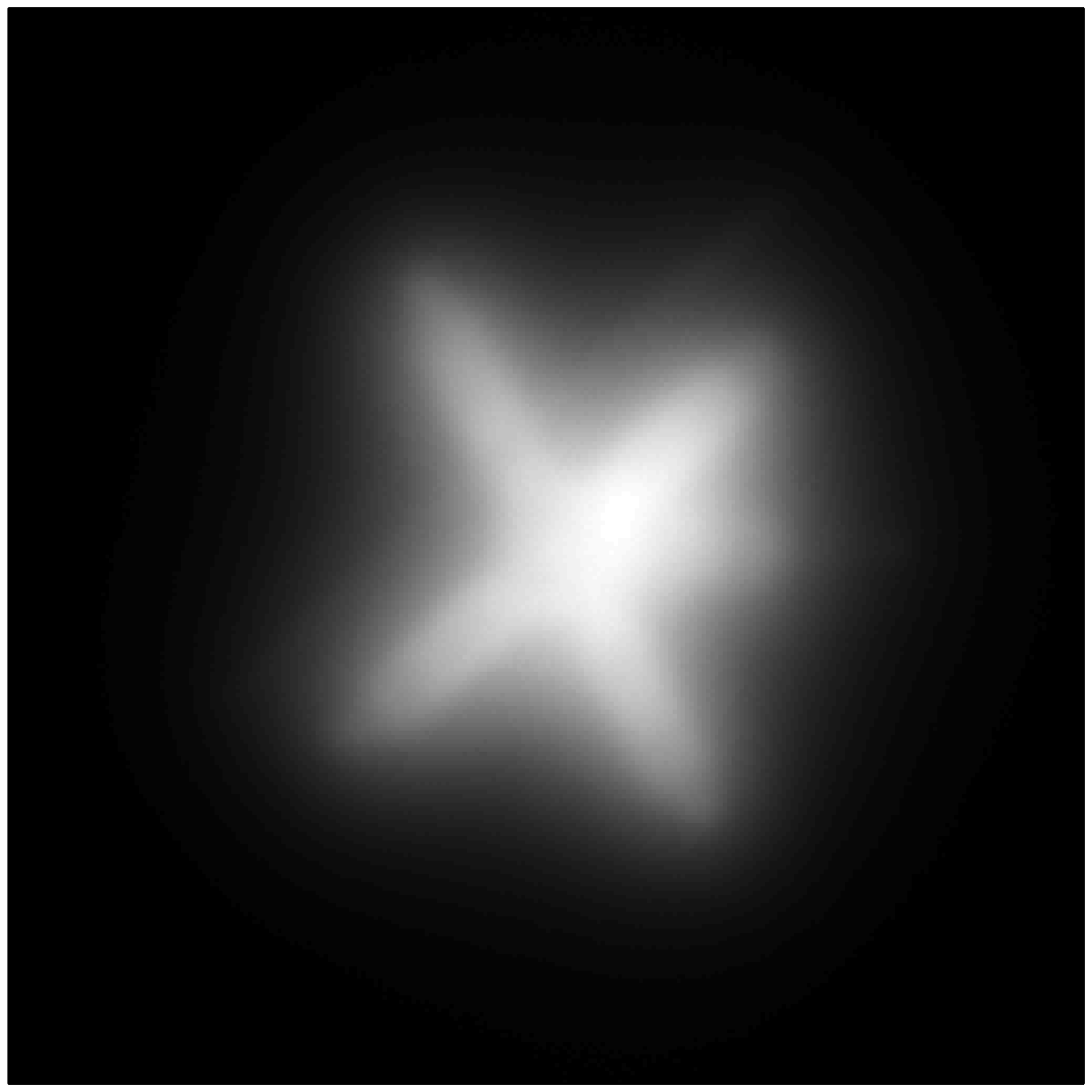}
\end{tabular}
\end{center}
\caption{Satellite test problems: original object (top left), S2550 (top right), S25500 (bottom left) and S255000 (bottom right) blurred and noisy images.}
\label{fig:3}
\end{figure}

The choice of the regularization parameter in Poisson data inversion is an active field and several different strategies have been proposed in the last years \cite{Bardsley2009,Bertero2010,Carlavan2011,Carlavan2012,Teuber2013}. Here we consider that proposed by Bertero et al. \cite{Bertero2010}, which consists of selecting the value of $\nu$ in \eqref{KL_R} such that 
\begin{equation}\label{discr_eq}
\mathcal{D}_A(\x_\nu) \equiv \frac{2}{n}KL(\x_\nu) = \eta,
\end{equation}
where $\eta$ is a given number close to 1 \cite{Bertero2010,Bonettini2014} (here we will assume $\eta=1$). In particular, in \cite{Zanni2015} the authors introduced an effective secant-type solver for the discrepancy equation \eqref{discr_eq}, called modified Dai-Fletcher (MDF) method, able to reduce the number of required solutions of problems \eqref{KL_R}. At each step of the secant method, an approximation of the solution of problem \eqref{KL_R} for a given value of $\nu$ is provided by running an optimization method until the stopping criterium
\begin{equation}\label{stcr}
|f(\x^{(k)})-f(\x^{(k-1)})| \leq \varepsilon |f(\x^{(k)})|,
\end{equation}
where $\varepsilon = 5 \times 10^{-8}$, is satisfied or when a maximum number of iterations equal to 5000 is reached. In this section we consider again the KL + HS model \eqref{KLHS} and we investigate the impact of (some of) the strategies used for the previous tests within this automatic scheme for the choice of $\nu$. In particular, we restrict our analysis to the GP, SGP$^*$ and SGP methods equipped with the Ritz-like steplengths, and the PidSplit+ algorithm with the adaptive choice of its parameter $\gamma$ described in \cite[equation (24)]{Zanni2015}, which resulted to be less dependent on the parameter settings than the standard approach.\\
The test problems we considered are based on the Satellite dataset already used in several papers and available at www.mathcs.emory.edu/$\sim$nagy/RestoreTools/index.html. The original image is sized $256 \times 256$ and assumes values in the range $[0,2550]$. The blurred image has been obtained by convolving the object with a point spread function simulating a ground-based telescope response, and a constant background $b=10$ has been added to the resulting image before introducing Poisson noise. Two further datasets have been obtained by multiplying object and background by factors of 10 and 100 before the blurring step. The three test sets will be denoted by S2550, S25500 and S255000 and the corrupted images are shown in figure \ref{fig:3} together with the original one. As concerns the parameter $\rho$ defining the HS regularization term, we followed the suggestion in \cite{Zanni2015} and set $\rho=10^{-4}\max(\ve g)$.\\

\noindent The results obtained by the algorithms are shown in table \ref{table}, where we reported the number of steps $k$ of the secant-based method required to satisfy either the relation
\begin{equation*}
|\mathcal{D}_A(\x_{\nu_k}) - \eta| \leq \varepsilon_1
\end{equation*}
or both the inequalities
\begin{equation*}
|\nu_k - \nu_{k-1}| \leq \varepsilon_2\nu_k \qquad ; \qquad |\mathcal{D}_A(\x_{\nu_k}) - \eta| \leq 10\varepsilon_1,
\end{equation*}
being $\varepsilon_1 = 5 \times 10^{-4}$ and $\varepsilon_2 = 5 \times 10^{-3}$, the total number of iterations $k_{\rm{tot}}$ performed by each method in the $k$ steps, the final regularization parameter $\nu_k$, the relative reconstruction error between $\x_{\nu_k}$ and $\xstar$ and the execution time in seconds. These numerical experiments has been carried out on a Dual CPU Intel(R) Xeon(R) X5690 at 3.47GHz with 188 GB RAM (see also fermi.unife.it) in a Matlab2013a environment.\\

\begin{table}[h]
\caption{\label{table}Results obtained with GP, SGP$^*$, SGP and PidSlit+ on the three Satellite datasets. Here $k$ denotes the number of steps of the secant-based method proposed in \cite{Zanni2015}, $k_{\rm{tot}}$ the total number of iterations performed by each method in the $k$ steps, $\nu_k$ the final regularization parameter, ``err'' the relative reconstruction error between $\x_{\nu_k}$ and $\xstar$ and ``time'' the execution time in seconds.}
\begin{indented}
\item[]\begin{tabular}{@{}ccccccc}
\br
Test problem             & Algorithm   & $k$  & $k_{\rm{tot}}$ & $\nu_k$ & err   & time       \\
\mr
\multirow{4}{*}{S2550}   & PidSplit+   & 8    & 1801       & 5.42e-05   & 0.322  & 67.84      \\
                         & GPRitz      & 16   & 3297       & 2.39e-04   & 0.308  & 56.55      \\
                         & SGPRitz$^*$ & 13   & 3091       & 6.29e-04   & 0.301  & 51.36      \\
                         & SGPRitz     & 7    & 2223       & 5.68e-04   & 0.305  & 36.52      \\
\mr
\multirow{4}{*}{S25500}  & PidSplit+   & 17   & 3486       & 1.00e-06   & 0.280  & 131.8      \\
                         & GPRitz      & 41   & 10388      & 1.00e-41   & 0.433  & 193.0      \\
                         & SGPRitz$^*$ & 15   & 7753       & 8.27e-05   & 0.257  & 130.3      \\
                         & SGPRitz     & 13   & 6189       & 9.91e-05   & 0.260  & 113.9      \\
\mr
\multirow{4}{*}{S255000} & PidSplit+   & 13   & 4108       & 7.97e-07   & 0.238  & 145.6      \\
                         & GPRitz      & 40   & 31150      & 1.00e-11   & 0.777  & 602.5      \\
                         & SGPRitz$^*$ & 16   & 10851      & 8.09e-06   & 0.230  & 189.7      \\
                         & SGPRitz     & 7    & 5341       & 1.00e-05   & 0.237  & 105.6      \\
\br											
\end{tabular}
\end{indented}
\end{table}

\noindent The performances summarized in table \ref{table} confirm what already observed in the previous section, since SGP equipped with the scaling matrices with variable bounds succeeds in reducing the overall number of iterations required to provide the regularization parameter and the corresponding reconstruction if compared with SGP with fixed bounds for the scaling matrices or GP (which, in two of the three tests, often fails in satisfying the stopping criterium \eqref{stcr} within the maximum number of iterations allowed). As concerns the comparison with PidSplit+, we can observe that the number of iterations performed by this latter strategy is lower than that of SGP, but the higher cost per iteration which characterizes PidSplit+ makes the procedure more expensive in terms of total CPU time with respect to the SGP method. 
\section{Conclusions}\label{sec:concl}
In this paper we revisited the SGP method, originally published in 2009 and exploited in the successive years in several inverse problems as image denoising/deblurring, Fourier-based image reconstruction, blind deconvolution, system identification and non-negative matrix factorization, with several applications in astronomy, microscopy and engineering. Despite all the good numerical results provided in solving these problems, the only theoretical convergence result proved so far is the stationarity of any limit point of the sequence generated by SGP. In this paper we showed that stronger results can be proved in the convex case, if the sequence of scaling matrices characterizing the SGP iterations is chosen as convergent to the identity matrix at a certain rate. Moreover, in the same setting we provided also a convergence rate estimate on the objective function values, as provided in the literature for several other optimization methods. Some numerical tests showed that the specific rule introduced on the scaling matrices to prove the theoretical convergence results helps also to improve the performances of the method, making SGP competitive also with methods for which the ${\mathcal O}(1/{k^2})$ convergence rate has been demonstrated.
\section*{Acknowledgments}
This work has been partially supported by MIUR (Italian Ministry for University and Research), under the projects FIRB - Futuro in Ricerca 2012, contract RBFR12M3AC, and PRIN 2012, contract 2012MTE38N. The Italian GNCS - INdAM (Gruppo Nazionale per il Calcolo Scientifico - Istituto Nazionale di Alta Matematica) is also acknowledged.
\section*{References}
\bibliographystyle{unsrt}
\bibliography{biblio_Marco}
\end{document}